\documentclass[12pt]{amsart}
\usepackage{amsfonts, amssymb,amsmath, amsthm, amsxtra, latexsym, mathrsfs, amscd}
\usepackage{amsmath, amssymb, amsthm, times, float}
\usepackage{graphics}
\usepackage[latin1]{inputenc} %

\newtheorem{theo+}{Theorem}[section]
\newtheorem{prop+}[theo+]{Proposition}
\newtheorem{coro+}[theo+]{Corollary}
\newtheorem{lemm+} [theo+]{Lemma}
\newtheorem{deep+}  [theo+]  {Deep Result}
\newtheorem{fact+}  [theo+]  {Fact}
\theoremstyle{definition}
\newtheorem{exam+}  [theo+]  {Example}
\newtheorem{rema+}  [theo+]  {Remark}
\newtheorem{defi+}  [theo+]  {Definition}
\newtheorem{xca+}[theo+]{Exercise}

\numberwithin{equation}{section}

\usepackage{color}

\def\draft{\centerline{(Draft {\the \day}/{\the\month} \the \year.)}}
\bigskip

\def\refn#1.#2{\expandafter\def\csname#1\endcsname{[#2]}}
\def\refnr#1.{\csname#1\endcsname}

\def\fgl{\mathfrak gl}

\def\fsl{\mathfrak sl}

\def\a{\alpha}

\def\Claminv2{|C(\Lambda)|^{-2}}

\def\de{d\varepsilon}

\def\Aa2D{A^{\a,2}(D)}
\def\bAa2D{\overline{A^{\a,2}(D)}}
\def\Ab2D{A^{\beta,2}(D)}
\def\bAb2D{\overline{A^{\beta,2}(D)}}

\def\Norm#1_#2{\Vert#1\Vert_{#2}}

\def\phipl12{\phi_{p_{l_1}, p_{l_2}}}
\def\phip01{\phi_{p_{0}, p_{0}}}

\def\a{\alpha}

\def\Claminv2{|C(\Lambda)|^{-2}}

\def\tr{\operatorname{tr}}

\def\End{\operatorname{End}}

\def\de{d\varepsilon}

\def\Aa2D{A^{\a,2}(D)}
\def\bAa2D{\overline{A^{\a,2}(D)}}
\def\Ab2D{A^{\beta,2}(D)}
\def\bAb2D{\overline{A^{\beta,2}(D)}}

\def\phipl12{\phi_{p_{l_1}, p_{l_2}}}
\def\phip01{\phi_{p_{0}, p_{0}}}

\def\br{\mathbb R}

\def\alg/{algebra}

\def\Alg/{Algebra}

\def\alt/{alternative} % \alt/ly
\def\anal/{analytic}
\def\analfunc/{\anal/\ \func/}
\def\Ans/{\it Answer. \normal}
\def\ass/{associative}
\def\nass/{non-\ass/}
\def\autom/{automorphism}
\def\homom/{homomorphism}
\def\isom/{isomorphism}
\def\bdd/{bounded}
\def\Bdd/{Bounded}
\def\bddsymdom/{bounded \sym/ \dom/}
\def\Cartdom/{Cartan \dom/}
\def\bdry/{boundary}
\def\bsd/{\bdd/ \symdom/}
\def\bv/{boundary value}
\def\cf/{{\it cf}\.}
\def\Cf/{{\it Cf}\.}
\def\charr/{character}
\def\coeff/{coefficient}
\def\comm/{commutative}
\def\cpct/{compact}
\def\compl/{complex}
\def\comp/{complex}
\def\Comp/{Complex}
\def\conf/{conformal}
\def\conj/{conjugate}
\def\conn/{connect}
\def\cont/{continuous}
\def\conv/{converge} % \conv/nce \conv/nt
\def\convc/{convergence}
\def\convt/{convergent}
\def\convx/{convex}
\def\coord/{coordinate}
\def\lcoord/{local coordinate}
\def\Corr/{Corresponding}
\def\corr/{corresponding}
\def\corrd/{correspond}
\def\cov/{covariant}
\def\decomp/{decomposition}
\def\deco/{decompose}
\def\diff/{different} % \diff/iable \diff/ial
\def\Diff/{Different} % \Diff/able \Diff/ial
\def\dimn/{dimension} % \dimen/al
\def\distr/{distribution} % \distr/al
\def\div/{diverge} % \div/nt
\def\dom/{domain}
\def\eg/{\hbox{\it e.g}\.}
\def\eigenf/{eigen\-\func/}
\def\eigensp/{eigen\-space}
\def\eigenv/{eigen\-value}
\def\eq/{equation}
\def\equa/{equation}
\def\de/{\diff/ial \equa/}
\def\do/{\diff/ial operator}
\def\ode/{ordinary \de/}
\def\pde/{partial \de/}
\def\pdo/{partial \diff/ial operator}
\def\psdo/{pseudo \diff/ial operator}
\def\fin/{finite}
\def\Ex/{\it Example.\ \normal}
\def\Exnr#1/{\it Example #1.\ \normal}
\def\foll/{follow}
\def\follg/{following}
\def\Follg/{Following}
\def\func/{function}
\def\Func/{Function}
\def\Fonc/{Fonc\-tion}
\def\fonc/{fonc\-tion}
\def\Funk/{Funk\-tion}
\def\funk/{Funk\-tion}
\def\gen/{general}
\def\har/{harmonic}
\def\Hint/{\it Hint. \normal}
\def\hist/{historic}
\def\histcl/{historical}
\def\hol/{holo\-morphic}
\def\homog/{ho\-mo\-ge\-ne\-ous}
\def\hyp/{hyper\-bolic}
\def\hyperg/{hyper\-geometric}
\def\ie/{\hbox{\it i.e.}}
\def\iff/{if and only if}
\def\ineq/{inequality}
\def\infra/{{\it inf\-ra}}
\def\ultra/{{\it ult\-ra}}
\def\Inpart/{In particular}
\def\inpart/{in particular}
\def\instof/{instead of}
\def\interps/{interpolation space}
\def\interp/{interpolation}
\def\Interp/{Interpolation}
\def\interpr/{Interpretation}
\def\Intr/{Introduction}
\def\intv/{interval}
\def\inv/{invariant}
\def\invc/{invariance}
\def\Iowords/{In other words}
\def\iowords/{in other words}
\def\ipr/{inner product}
\def\irred/{irreducible}
\def\lb/{line bundle}
\def\lin/{linear}
\def\lhs/{left hand side}
\def\rhs/{right hand side}
\def\loc/{local}
\def\math/{mathematic}
\def\mathcn/{\math/ian}
\def\manif/{manifold}
\def\meas/{measure}
\def\measl/{measurable}
\def\mero/{mero\-morphic}
\def\mon/{monomial}
\def\monog/{monogenic}
\def\mult/{multiple}
\def\multy/{multiply}
\def\multn/{multiplication}
\def\nas/{necessary and sufficient}
\def\nbd/{neighborhood}
\def\neg/{negative}
\def\nondeg/{nondegenerate}
\def\Oohand/{On the other hand}
\def\oohand/{on the other hand}
\def\Oonhand/{On the one hand}
\def\oonhand/{on the one hand}
\def\oper/{operator}
\def\orth/{ortho\-gonal}
\def\orthon/{ortho\-normal}
\def\otoh/{on the other hand}
\def\quat/{quaternion}
\def\pp/{\hbox{a. e.}}
\def\psh/{plurisubharmonic}
\def\pol/{polynomial}
\def\pot/{potential}
\def\pos/{positive}
\def\princ/{principle}
\def\prob/{probability}
\def\proj/{projective}
\def\projn/{projection}
\def\Proof/{\it Proof:\normal}
\def\Rem/{\it Remark\normal}
\def\Remnr#1/{\it Remark\ \normal #1. }
\def\rep/{representation}
\def\reps/{representations}
\def\meta/{metaplectic representation}
\def\repr/{reproducing}
\def\reprker/{reproducing kernel}
\def\resp/{respective} % \resp/ly
\def\resply/{respectively}
\def\restr/{restriction}
\def\sa/{self-adjoint}
\def\st/{such that}

\def\sol/{solution}
\def\ru/{space}
\def\sph/{spherical}
\def\ssp/{sub\ru/}
\def\sym/{symmetric}
\def\Sym/{Symmetric}
\def\symb/{symbol}
\def\symbc/{symbolic}
\def\symdom/{\sym/ domain}
\def\symp/{symplectic}
\def\Theor#1/{\fet Theorem #1.\ \normal}
\def\Lem#1/{\fet Lemma #1.\ \normal}
\def\Lemma/{\fet Lemma.\ \normal}
\def\topl/{topology}
\def\topll/{topological}
\def\transf/{transform}
\def\transl/{translation}
\def\transfn/{transformation}
\def\transv/{transvectant}
\def\trig/{trigonometric}
\def\tril/{trilinear}
\def\trilf/{trilinear form}
\def\uhp/{upper halfplane}
\def\uhs/{upper halfspace}
\def\vb/{vector bundle}
\def\vf/{vector field}
\def\vsp/{vector space}
\def\wrt/{with respect to}
\def\Wlog/{Without loss of generality}
\def\a{\alpha}

\def\Ab/{Abel}
\def\Ban/{Banach}
\def\Bansp/{\Ban/ space}
\def\Belt/{Bel\-tra\-mi}
\def\Berg/{Berg\-man}
\def\Bern/{Ber\-nou\-lli}
\def\Berz/{Berezin}
\def\Bess/{Bessel}
\def\Cart/{Car\-tan}
\def\Cay/{Cay\-ley}
\def\CG/{Clebsch-Gordan}
\def\Cl/{Clifford}
\def\CR/{Cauchy-Rie\-mann}
\def\Dir/{Dirichlet}
\def\Eucl/{Euclide}
\def\Eucln/{Euclidean}
\def\F/{Fourier}
\def\Hank/{Hankel}
\def\Hankf/{\Hank/ form}
\def\Herm/{Hermite}
\def\Hilb/{Hilbert}
\def\Hilbs/{Hilbert space}
\def\Hilbsp/{Hilbert space}
\def\HS/{Hilbert-Schmidt}
\def\Lag/{La\-grange}
\def\Lap/{La\-place}
\def\LapBelt/{\Lap/-\Belt/}
\def\Leb/{Lebesgue}
\def\Marc/{Mar\-cin\-kie\-wicz}
\def\Moeb/{Moebius}
\def\Moebt/{Moebius transformation}
\def\Moebtransfn/{Moebius transformation}
\def\Pla/{Plan\-che\-rel}
\def\Poin/{Poin\-car\'e}
\def\Riem/{Rie\-mann}
\def\Riemn/{\Riem/ian}
\def\psRiemn/{pseudo-\Riem/ian}
\def\Riems/{Rie\-mann surface}
\def\Schroe/{Schr\"odinger}
\def\Weier/{Weier\-strass}
\def\anal/{analytic}
\def\bsd/{bounded symmetric domain  }
\def\bdd/{bounded}
\def\calc/{calculation}\def\conj{conjugate}
\def\calci/{calculating}\def\eg{e.g.}
\def\conj/{conjugate}
\def\deco/{decomposition}
\def\eg/{e.g.}
\def\fct/{function}
\def\gp/{group}
\def\hw/{highest weight}
\def\hwv/{highest weight vector}
\def\hwvs/{highest weight vectors}
\def\lw/{lowest weight}
\def\lwv/{lowest weight vector}
\def\lwvs/{lowest weight vectors}
\def\hds/{holomorphic discrete series}
\def\iff/{if and only if}
\def\inv/{invariant}
\def\irrde/{irreducible decomposition}
\def\meas/{measure}
\def\transf/{transform}
\def\rep/{representation}
\def\resp/{respectively}
\def\inters/{intertwines}
\def\interg/{intertwining}
\def\meta/{metaplectic representation}
\def\qu/{quaternion}
\def\rep/{representation}
\def\symdom/{ symmetric domain}
\def\st/{such that}
\def\shd/{subhead}
\def\transf/{transform}
\def\wrt/{with respect to}

\def\ra{\rightarrow}

\def\Norm#1#2#3{\Vert#1\Vert^{#3}_{{#2}}}

\def\tr{\operatorname{tr}}

\begin{document}

\def\abstractname{Abstract}
\def\chrefname{References}
%\today

\title[Weil's local rigidity theorem
]{Convex real projective structures and Weil's local rigidity Theorem
}
\author{Inkang Kim and  Genkai Zhang}
%%%%

\address{School of Mathematics,
KIAS, Heogiro 85, Dongdaemun-gu
Seoul, 130-722, Republic of Korea.
\text{Email: inkang@kias.re.kr}
}
\address{Mathematical Sciences, Chalmers University of Technology and
Mathematical Sciences, G\"oteborg University, SE-412 96 G\"oteborg, Sweden.
\text{Email: genkai@chalmers.se}
}
\footnotetext[3]{Research partially supported by
STINT-NRF grant (2011-0031291). Research by G. Zhang is supported partially
 by the Swedish
Science Council (VR). I. Kim gratefully acknowledges the partial support
of grant  (NRF-2014R1A2A2A01005574) and a warm support of
Chalmers University of Technology during his stay.
}
%\draft
\begin{abstract} For an $n$-dimensional real hyperbolic manifold $M$,
we  calculate the Zariski tangent space of a character variety $\chi(\pi_1(M),SL(n+1,\mathbb R)), n>2$ at Fuchisan loci to show that the tangent space consists of cubic forms. Furthermore we prove the Weil's local rigidity theorem for uniforml hyperbolic lattices using real projective structures.
\end{abstract}

\maketitle

\baselineskip 1.35pc

\section{Introduction}

A flat projective structure on an $n$-dimensional manifold $M$ is a $(\mathbb{RP}^n, PSL(n+1,\br))$-structure, i.e., there exists a maximal atlas on $M$ whose transition maps are restrictions to open sets in $\mathbb{RP}^n$ of elements in $PSL(n+1,\br)$. Then there exist a natural holonomy map $\rho:\pi_1(M)\rightarrow PSL(n+1,\br)$ and a developing map $f:\tilde M\rightarrow \mathbb{RP}^n$ such that
$$\forall x\in\tilde M, \forall \gamma\in\pi_1(M), \ f(\gamma x)=\rho(\gamma)f(x).$$

In this paper, since we will consider projective structures deformed from  hyperbolic structures, all holonomy representations will lift to $SL(n+1,\br)$.
An $\mathbb{RP}^n$-structure is convex if the developing map is a homeomorphism onto a convex domain in $\mathbb{RP}^n$.  It is {\it properly convex} if the domain is included in a compact convex set of an affine chart, {\it strictly convex} if the convex set is strictly convex.
Surprisingly, while many people are working on global structures of the Hitchin component, it seems that the local structure
has been neglected. This is the starting point of this article.
We shall first compute  the cohomology
$H^1(\pi_1(M), \rho, \fsl(n+1, \mathbb R))$
of a Fuchsian point $\rho$ which  corresponds to a hyperbolic structure $\pi_1(M)\ra SO(n,1)\subset SL(n+1,\br)$. See  Section \ref{projectivestructure} for details. The cohomology
is described in terms of quadratic and cubic
forms.  We shall use
the result of Labourie \cite{La} where he
proved that a convex projective
flat structure on $M$ defines
a Riemannian metric and a cubic form on $M$.
%As a consequence of our result for $n=2$  we prove
%a Eichler-Shimura isomorphism
%characterizing the tangent space
%$H^1(\pi_1(M), \rho, \fsl(3, \mathbb R))$
%of Hitchin component
%in terms of quadratic and cubic holomorphic forms.
\begin{theo+}
Let $\rho: \pi_1(M)\to SO(n,1)\subset SL(n+1, \mathbb R), n>2,$
be a representation defining
a real hyperbolic   structure
on a closed n-manifold $M$.
Let
$\alpha\in H^1(\pi_1(M), \rho, \fsl(n+1, \mathbb R))$. Then
 $\alpha$ is represented by a cubic form.
\end{theo+}
For $n=2$,  an element in $H^1(\pi_1(M),\rho, \fsl(3,\mathbb R))$ is represented by a sum of a quadratic form and a cubic form when
$\rho$ defines a convex real projective structure. 
This is due to \cite{Loftin, La}. In this case both
the global and local structures have
been studied intensively. Recently  \cite{KZ-hitchin} 
we have been able 
to construct mapping class group invariant K\"ahler metric
on the Hitchin component of $SL(3, \mathbb R)$, 
this is also part of our motivation of the present paper.
We also mention that Labourie \cite{La} has computed
the cohomology $H^1(\pi_1(M), \mathbb R^{3})$
where $\pi_1(M)$ acts on  $\mathbb R^{3}$
through $\rho$ and the defining representation of $SL(3, \mathbb R)$.

As a corollary of our technique, we show the Weil's local rigidity theorem for uniform real hyperbolic lattices for dimension $>2$.

\begin{theo+}\label{Weil-thm} Let $M=\Gamma\backslash SO^0(n, 1)/SO(n)$
be a compact hyperbolic manifold.
If $n>2$ then   $ H^1(\Gamma, \mathfrak{so}(n, 1))=0$.
\end{theo+}

{\bf Acknowledgement}
 We are grateful
for the anonymous referee for  the careful reading of an earlier version
of this paper and many useful comments.

\section{Tangent space at Fuchsian locus
of convex projective structures on a manifold $M$
}
\subsection{Projective structure}\label{projectivestructure}

The notion of projective structures can be formulated in terms of a flat connection as follows, see \cite{La} for details. Consider a trivial bundle $E=M\times \br^{n+1}$ where $M$ is an $n$-dimensional manifold. Let $\omega$ be a volume form on $\br^{n+1}$ and let $\nabla$ be a flat connection on $E$ preserving $\omega$.  Let $\rho$ be the holonomy  representation of $\nabla$. A section $u$ of $E$ is identified with a $\rho$-equivariant map from $\tilde M$ to $\br^{n+1}$. A section $u$ is said to be $\nabla$-immersed if the $n$-form $\Omega_u$ defined by
$$\Omega_u(X_1,\cdots,X_n)=\omega(\nabla_{X_1} u,\cdots,\nabla_{X_n} u, u)$$ is non-degenerate. Then $u$ is $\nabla$-immersed if it is a non-vanishing section and if the associated $\rho$-equivariant map $[u]$ from $\tilde M$ to $\mathbb{RP}^n$ is an immersion.

Hence it follows that such a pair $(\nabla, u)$ gives rise to a flat projective structure.
Labourie reformulated this  as a pair of torsion free connection $\nabla^T$ on $M$ with a symmetric 2-tensor $h$ on $M$ as follows. One can associate a connection $\nabla$ on $E=TM\oplus \mathcal L$, where $\mathcal L$ is a trivial bundle $M\times \br$, defined by
\begin{equation}
 \label{nabla-E1}
\nabla_X
\begin{pmatrix} Y\\
\lambda
\end{pmatrix}
=
\begin{pmatrix} \nabla_X^T
 & X
\\
h(X, \cdot) & L_X
\end{pmatrix}
\begin{pmatrix} Y\\
\lambda
\end{pmatrix}
=
\begin{pmatrix}
 \nabla_X^T Y +\lambda X
\\
h(X, Y)+  L_X(\lambda)
\end{pmatrix}.
\end{equation}
Here $L_X(\lambda)=X\lambda$ denotes the differentiation.
Labourie \cite{La} showed that if $\nabla$ is flat and $\nabla^T$ preserves the volume form defined by $h$, then $\nabla$ gives rise to a flat projective structure. He further showed that $h$ is positive definite if the structure is properly convex. We will use this final form of projective structure in this paper to
carry out the explicit calculations.
\subsection{Tangent space of convex projective structures}
Let $M$ be an $n$-dimensional manifold and $\Gamma=\pi_1(M)$
its fundamental group. Let $\rho$ be a representation
of $\Gamma$ into $SL(n+1, \mathbb R)$ defining
a  convex projective structure on $M$.
There is   a flat connection $
\nabla$ on the  associated
$\mathbb R^{n+1}$-bundle $E$
preserving a  volume form as in the previous section.

The flat connection $\nabla$ on $ E$
defines also a connection
 on the dual bundle
$ E^\ast$,
$\fgl(E)= E\otimes E^\ast
$, the bundle of endomorphisms
of $E$, and further on $\fsl(E)$,
 the trace free endomorphisms,
since $\nabla$ preserve the volume form on $ E$,
by the Leibniz rule
and the commutative relation with the contractions.

Write temporarily $\mathfrak F$
for any of these flat bundles with fiber space $F$.
We fix convention that if $\mathfrak F= E^\ast$
or
$\fsl( E)
\subset
\fgl( E)= E\otimes  E^\ast
$, we write then  a $F$-valued
one-form $\alpha$ as $\alpha: (X, y)\to \alpha(X)y$ with the first argument
being tangent vector and the second argument being element of $F$.
 For conceptual
clarity we recall that
given $\nabla^{E}$ on $E$
the connection
$\nabla^{E^\ast}$
on $E^\ast$-valued sections is
defined by
the equation
$$
X(\alpha(y))=
(\nabla^{E^\ast}_X
\alpha)(y)+
\alpha
(\nabla_X^{ E} y);
$$
whereas the connection on
sections of $\fsl( E)$ is defined by
\begin{equation}
\label{nabla}
\nabla^E_X(\alpha (y))=
(\nabla_X^{sl} \alpha)(y)
 +
\alpha(\nabla_X^{ E} y).
\end{equation}
We shall abbreviate them all as $
\nabla_X$.

The flat connection $\nabla$
on $E$ as well as its induced
connection $\nabla^T$
induces  exterior differentiation
 $d^{\nabla}
$ and  $d^{\nabla^T}
$ on
 $1$-forms, locally defined as
$$d^\nabla (\sum  \omega_i dx_i)=\sum \nabla\omega_i dx_i,$$ where $\omega_i$ is a local section.
  For notational convenience
we shall write all of them just as $d^\nabla$;
 no confusion
would arise as it will be clear from the context  which
sections they are acting on.
We will freely write a $(0, 2)$-tensor
as $\alpha(X)Y=
\alpha(X, Y)$.
We shall need the following
formula for the exterior differentiation
 on a $End(E)$-valued
one-form:
\begin{equation}
\label{ext-diff}
\begin{split}
&\qquad (d^{\nabla}\alpha)(X, Z)y\\
&=({\nabla}_X \alpha)(Z)y
-({\nabla}_Z \alpha)(X)y
\\
&={\nabla}_X (\alpha(Z)y) -\alpha(\nabla_X Z)y
-\alpha(Z) (\nabla_X
y) \\
&\qquad -\left({\nabla}_Z (\alpha(X)y) -\alpha(\nabla_Z X)y
-\alpha(X) (\nabla_Z y)\right)
\\
&=
{\nabla}_X (\alpha(Z)y) -\alpha(Z) (\nabla_X
y) -({\nabla}_Z (\alpha(X)y) -\alpha(X) (\nabla_Z y))
-\alpha([X, Z])y
\end{split}
\end{equation}
since $\nabla$ is flat, in particular $\nabla^T$ is torsion free.

We shall describe the cohomology in terms of some symmetry
conditions of certain tensors. Let $g$
be a Riemannian metric on $M$ with
the Levi-Civita connection $\nabla^g$, the corresponding exterior
differentiation being denoted by $d^g$.
We consider following
three conditions for  $(0, 2)$-tensors $\alpha$, in which case $\alpha$ is a quadratic form:
\begin{description}
\item[(q1)] $\alpha$ is symmetric,

\item[(q2)] trace-free with respect to $g$,

\item[(q3)] $ \alpha$ is $d^{\nabla^g}$ closed, $ d^{\nabla^g}\alpha=0, $
\end{description}

and the following four conditions
for a $End(TM)$-valued one-form, in which case $\alpha$ is a cubic form:

\begin{description}
\item[(c1)]
 $\alpha(X)Y
$ is symmetric in $X, Y$,
$ \alpha(X)Y =
\alpha(Y)X$,

\item[(c2)]
$\alpha(X)$ is symmetric with respect to $g$,
 $\alpha(X)^\ast=\alpha(X)$,

\item[(c3)]
$\alpha $ is $d^{\nabla^g}$
closed, $d^{\nabla^g}\alpha=0$, equivalently the cubic form
$g(\alpha(X)Y, Z)$ is closed,

\item[(c4)] $\alpha(X)$ is trace-free.

\end{description}

The following theorem is proved in \cite[Theorems 3.2.1 \& Proposition 4.2.3]{La}.

\begin{theo+}\label{Labourie}
Suppose $M$ admits a properly convex projective structure.
Then there is splitting
of the bundle
\begin{equation}
\label{split}
 E=TM\oplus \mathcal L
\end{equation} where $\mathcal L$ is a trivial bundle $M\times\br$
and a  Riemannian metric
$g$ on $M$ such that
 the flat
connection  $\nabla  $ on $E$ is given
by
\begin{equation}
 \label{nabla-E}
\nabla_X
\begin{pmatrix} Y\\
\lambda
\end{pmatrix}
=
\begin{pmatrix} \nabla_X^T
 & X
\\
g(X, \cdot) & L_X
\end{pmatrix}
\begin{pmatrix} Y\\
\lambda
\end{pmatrix}
=
\begin{pmatrix}
 \nabla_X^T Y +\lambda X
\\
g(X, Y)+  L_X(\lambda)
\end{pmatrix}.
\end{equation}
Here  $L_X(\lambda)=X\lambda$  is the differentiation, $\nabla^T$ is a torsion-free connection on $TM$
preserving the volume form of the Riemannian metric
$g$ and $\nabla^g$ is the Levi-civita connection of $g$, such that the tensor $c$, defined by
$g((\nabla_X^T-\nabla^g_{ X})Y, Z)
=c(X, Y, Z)$,
satisfies the condition (c1-c4).
\end{theo+}

%More precisely the second part
%of the statement
%that
%$g((\nabla_X^T-\nabla^g_{ X})Y, Z)
%=c(X, Y, Z)$
%satisfies the condition (c1-c4) follows
%from the proof of Proposition 4.2.3 which
%is done for surfaces there. However it depends
%only the flat condition of the connection
%and the curvature computation along with
%the given Riemannian metric.

In other words such a flat connection
$\nabla$
determines a  metric $g$
on $TM$ and thus on  the bundle $TM\oplus \mathcal L$,
$$
g\oplus |\cdot|^2:
(X, \lambda)\mapsto g(X, X)+\lambda^2
$$
such that the
connection
$\nabla$
takes the form
$$
\nabla_X=
\begin{bmatrix}
\nabla^g_X &0
\\
0&L_X
\end{bmatrix}
+
\begin{bmatrix}
Q(X)&X
\\
X_g^\sharp &0
\end{bmatrix}
$$
where $\nabla^T_X-\nabla^g_X=Q(X)$, hence the first part is a skew-symmetric (i.e. orthogonal) connection
and the second part is a symmetric form.

We shall now compute the tangent space
of the deformation space
at  Fuchsian locus induced
from the natural inclusion $SO(n, 1)\subset SL(n+1, \mathbb R)$.
A hyperbolic $n$-manifold can be viewed as a real projective manifold and it can be deformed
inside  the convex real projective structures for any dimension $n\geq 2$. The component containing
the hyperbolic structures constitutes the strictly convex real projective structures \cite{Be,CG,Kim}. Indeed, any convex real projective structures whose holonomy group is not contained in $SO(n,1)$ has a Zariski dense holonomy group in $SL(n+1,\mathbb R)$, see \cite{Be1}.

 Hence from now on, we shall assume $Q=0$ and
then $g$ is a hyperbolic metric on $M$ with constant
curvature $-1$.
All the covariant differentiations
below
will be the
 one induced by $\nabla^g$.

\begin{theo+}\label{2.2}
Let $\rho: \pi_1(M)\to SO(n,1)\subset SL(n+1, \mathbb R),\ n>2$
be a representation defining
a hyperbolic   structure
on the compact n-manifold $M$.
Let $g$ be the hyperbolic metric
determined by $\nabla$.
Then there exists an injective
map  from $ H^1(\pi_1(M), \rho, \fsl(n+1, \mathbb R))$
into the space of cubic forms satisfying (c1-c4).
\end{theo+}

The proof will be divided into several steps.

Let $\alpha$ be a one-form representing
an element
of
$ H^1(\pi_1(M), \rho, \fsl(n+1, \mathbb R))$, realized as
an element of $\Omega^1(M, \fsl(E))$. Write it
as
$$
X\mapsto \alpha(X)=
\begin{bmatrix}
A(X)&B(X)
\\
C(X)&D(X)
\end{bmatrix}\in \fsl(E),
$$
under the splitting (\ref{split}),
where $A\in \Omega^1(M,\End(TM))$,
$B\in \Omega^1(M, TM)$,
 $C\in C^\infty(M, T^*M\otimes T^*M)$ and $X\to D(X)=-\tr A(X)$
is a one-form.

We first observe for
the covariant differentiation
of a section of $\fsl( E)$
$$
u=
\begin{pmatrix}
a&b
\\
c&e
\end{pmatrix}
$$
is the
one-form
\begin{equation}
\label{nabla-u}
X\to \nabla_X u
=\begin{pmatrix}
(\nabla_X^g a)  \cdot  + c(\cdot) X - g(X, \cdot )b
  & \nabla_X^g b +e X - a(X)
\\
g(X, a\cdot) + (\nabla_X^g c)(\cdot) -e(X, \cdot)
& (Xe) +g(X, b) -c(X)
\end{pmatrix},
\end{equation}
acting on a section $y=(Y, \lambda)$, where the dots denote the variable $Y$.
These are  exact forms.

\begin{lemm+}
\label{lemm23}
Up to  exact forms we can assume
$$
B =0, D=0.
$$
\end{lemm+}

\begin{proof} We prove first  that
we can  choose  $u$
so that $\alpha +\nabla u$
has its entries $B(X)=bX, D=0$ where $b$ is a scalar function.
Indeed we choose
$$
u=
\begin{pmatrix}
a&0
\\
c&e
\end{pmatrix},
$$
$a=B$, $c$ the one-form $c= de + D$. The
form $\alpha +\nabla u$ then
has the desired form:
\begin{equation*}
\begin{split}
X\mapsto \alpha(X)
+\nabla_X u
&=\begin{bmatrix}
A(X)
+
(\nabla_X^g B)\cdot  + c(\cdot ) X
&B(X) + eX- B(X)
\\
C(X)
+
g(X, a\cdot) +
(\nabla^g_X c)\cdot  - g(X, \cdot) e
&D(X)+
L_X (e) -c(X)
\end{bmatrix}
\\
&=\begin{bmatrix}
A'(X)& e X
\\
C'(X)&0
\end{bmatrix}
\end{split},
\end{equation*}
as claimed. We write the new form $\alpha +\nabla u$
as $\alpha$ with its entries $B=b \text{Id}$, $D=0$.

Next we take
$$
v=
\begin{pmatrix}
a_1\text {Id} &0
\\
c_1& -n a_1
\end{pmatrix}, \quad a_1:=\frac 1{n+1} b, \,\, c_1:=-n da_1=-\frac n{n+1} db,
$$
the same calculation above shows that
$\alpha +\nabla v$ has its $B=0$, $D=0$.
\end{proof}

We shall need the precise formula of the diagonal part in the lower
triangular form, keeping tract of the computations we find
$\alpha +\nabla (u+v)
$ has the form
\begin{equation}
\label{newA}
\begin{pmatrix} A(X) +\nabla_X^g B  +
(Xb) I +f X & 0\\
\ast & 0
\end{pmatrix}, \quad b=-\frac{1}{n+1}\tr B, \quad f= db
\end{equation}

\begin{lemm+} Let $\alpha$ be  lower triangular with
$B=0$, $D=0$. The covariant derivatives
$d^{g} A$ and
$d^{g} C$ are related to $A$ and $C$ by
\begin{equation}
\label{cnd-1}
\boxed{
0 =(d^{g} A)
(X, Z)(Y)
+(C(Z)Y)X-(C(X)Y)Z}
\end{equation}
and
\begin{equation}
\label{cnd-2}
\boxed{
0= g(X, A(Z)Y)-
g(Z, A(X)Y)+ (d^g C)(X, Z)(Y),}
\end{equation}
and there hold the symmetric relations:
\begin{equation}
\label{A-b-sym}
\boxed{A(X)Z =A(Z)X},
\end{equation}
\begin{equation}
\label{C-sym}
\boxed{C(X)Z =C(Z)X}.
\end{equation}
\end{lemm+}

\begin{proof}
We write an arbitrary section $y$ of $ E$
as $y=(Y, \lambda)=(y_T, y_n)$, the tangential
and respectively normal component.
To write down the condition on the  closedness  $d^\nabla\alpha(X,
Z)=0$ in terms of the components $A, C$,
we recall (\ref{ext-diff}).
The condition  $d^\nabla\alpha(X, Z)=0$
is then
\begin{equation}
\label{da=0}
\begin{split}
0&=d^\nabla\alpha(X, Z)y\\
&=v(X, Z; y) -v(Z, X; y) -\alpha([X, Z])y\\
&=
v(X, Z; y) -v(Z, X; y) -  \begin{pmatrix}
A([X, Z])Y
\\ C([X, Z])Y
  \end{pmatrix}.
\end{split}
\end{equation}
where
\begin{equation*}
\begin{split}
&\qquad v(X, Z; y):=\nabla_X(\alpha(Z)y)-\alpha(Z)(\nabla_X y)=
\\
&
\begin{pmatrix}
    \nabla_X^g(A(Z)Y ) + (C(Z)Y) X
-A(Z)(\nabla_X^g Y +\lambda X)\\
    g(X, A(Z)Y ) + X(C(Z)Y) - C(Z)(\nabla^g_X Y +\lambda X)
  \end{pmatrix}
\end{split}.
\end{equation*}
Here we have used the fact that $B=0, D=0$ in the computations.

Recall that $d^{g} A=
d^{\nabla^g} A
$ is
defined
by
\begin{equation}
\label{dA}
(d^{g} A)(X, Z)(Y)=$$$$
\nabla^g_X(A(Z)Y)-
A(Z)(\nabla^g_X Y)
-\left(\nabla^g_Z(A(X)Y)
 -A(X)(\nabla^g_Z Y)\right)
-A([X, Z])Y
\end{equation}
which is a well-defined $End(TM)$-valued
 2-form,
and
\begin{equation}\label{dc}
\begin{split}
&\phantom{=}(d^{g} C)(X, Z)(Y)
\\
&=
L_X(C(Z)Y) -
 C(Z)(\nabla_X^gY)
-\left(L_Z(C(X)Y) - C(X)(\nabla_Z^gY)\right)\\
&- C([X, Z])Y
\end{split}
\end{equation}
 is the Riemannian exterior
differential of the form $C$,
and is an element of $\Omega^2\otimes \Omega$. The first two equations
(\ref{cnd-1})
-(\ref{cnd-2})
 are obtained
from (\ref{da=0}) by
putting $y=(Y,0)$.

Correspondingly we have, taking
 $y=(0, 1)$,
\begin{equation}
\label{cnd-3}
A(X)Z
-A(Z)X
=0
\end{equation}
and
\begin{equation}
\label{cnd-4}C(X)Z -C(Z)X=0,
\end{equation}
resulting the  symmetric relations (\ref{A-b-sym})-
(\ref{C-sym}).
\end{proof}

Let $\alpha_0$ be the bilinear form
$$
\alpha_0(X, W)=C(X, W):=C(X)W
$$
and
$\alpha_1$ the $End(TM)$-valued one-form
$$
\alpha_1(Y)X:=  \frac 12(A(Y)X +A_g^*(Y)X).
$$
Then $\alpha_0(X,W)$ is symmetric in $X$ and $W$, hence it satisfies
(q1). $\alpha_1(Y)$ is symmetric with respect to $g$
and trace free, since $\alpha(Y)\in \mathfrak{sl}(E)$,
 $0=\tr \alpha(Y)=\tr A(Y) + D(Y)= \tr A(Y)$, hence it satisfies the conditions (c1) and (c4).

\begin{lemm+}\label{lemma2.5}
 Let $\alpha$ be of the above form with $B=0, D=0$.
Then we have $C=0$ and $\alpha_0=0$ and $\alpha_1$
satisfies the conditions (c1)-(c4) for $n>2$.
\end{lemm+}
\begin{proof}
%It follows from the definition that $\tr A(X)=0$, namely $\alpha_1$ satisfies  $(c4)$.
  The equation
(\ref{cnd-2}) combined with
(\ref{A-b-sym}) implies that
\begin{equation*}
 g(X, A(Y)Z)-g(Z, A(Y)X) +(d^g C)(X, Z)Y
=0.
\end{equation*}
In other words
\begin{equation}
\label{a-star-a}
(A_g^\ast(Y)-A(Y))X
=-((d^g C)(X, \cdot)(Y))_g^{\sharp}
\end{equation}
where the lowering of the index in the right hand side
is with respect to the second variable. Since $g$ is parallell
with respect to $\nabla^g$ then 
$((d^g C)(X, \cdot)(Y))_g^{\sharp}$
 is an exact $End(TM)$-valued
one form.
This is not obvious and requires proof. Indeed,
let  $C^\flat$ be the 
$End(TM)$-valued 0-form,
$$
 C^\flat
 (X)=\sum_i C(X, Z_i)Z_i
$$
where $\{Z_i\}$ is a local orthonormal frame of $TM$.
We claim that 
\begin{equation}
\label{d-comm-w-flat}
((d^g C)(X, \cdot)(Y))_g^{\sharp}=(d_Y^g  C^\flat)(X).
\end{equation}
By definition we have the identity section $\text{Id} =\sum_i Z_i\otimes Z_i^{\sharp}
$ and, 
$$
0=\nabla_Y^g \text{Id} =
\sum_i (\nabla_Y^g Z_i\otimes Z_i^{\sharp} +
Z_i\otimes \nabla_Y^g (Z_i^{\sharp})
=\sum_i 
\nabla_Y^g Z_i\otimes Z_i^{\sharp} +
\sum_i 
Z_i\otimes (\nabla_Y^g Z_i)^{\sharp},
$$
namely, for any $Z$
\begin{equation}
  \label{eq:ident-zero}
0=\sum_i 
 {g}\left(Z_i, Z\right) \nabla_Y^g Z_i
+ \sum_i    {g}\left(\nabla_Y^g Z_i, Z\right)Z_i.  
\end{equation}
Here we have used the fact that $\nabla^g_Y$ commutes with $\sharp$,
$\nabla_Y^g (Z_i^{\sharp})=(\nabla_Y^g Z_i)^{\sharp}$.
By definition, LHS of (\ref{d-comm-w-flat}) is
\begin{equation*}
  \begin{split}
LHS&=\sum_i((d^g C)(X, Z_i)(Y)) Z_i=
\sum_i((\nabla_Y^g C)(X, Z_i)) Z_i\\
& =\sum_i  Y(C(X, Z_i)) Z_i - \sum_i C(\nabla^g_YX,  Z_i)Z_i - \sum_i C(X, \nabla^g_YZ_i)Z_i.
  \end{split}
\end{equation*}
Here $C\in C^\infty(M, T^*M\otimes T^*M)$ is a zero form, hence $(d^g C)(X, Z_i)(Y)=(\nabla_Y^g C)(X, Z_i)$.
On the other hand,
\begin{equation*}
  \begin{split}
RHS&=(\nabla^g_Y  C^\flat)
(X)=
\nabla^g_Y (C^\flat(X))-  C^\flat (\nabla^g_Y X)
\\
&=\nabla^g_Y \left(\sum_i C(X, Z_i)Z_i\right)- 
\sum_i C (\nabla^g_Y X, Z_i)Z_i\\
&=\sum_i Y(C(X, Z_i))  Z_i +
\sum_i C(X, Z_i)  \nabla^g_YZ_i
-\sum_i C (\nabla^g_Y X, Z_i)Z_i.
  \end{split}
\end{equation*}
To treat the second term we compute its inner product with any $Z$; it is
\begin{equation*}
  \begin{split}
\sum_i C(X, Z_i)  g(\nabla^g_Y Z_i, Z)
&= C(X, \sum_i g(\nabla^g_Y Z_i, Z) Z_i)\\
&=- C(X, \sum_i g(Z_i, Z) \nabla_Y^gZ_i
)\\
&=-  g(\sum_i C(X,  \nabla_Y^gZ_i)Z_i, Z)
  \end{split}
\end{equation*}
where the second equality is by
(\ref{eq:ident-zero}). Hence 
$$
\sum_i (C(X, Z_i))  \nabla^g_Y Z_i=
- \sum_i C(X, \nabla^g_YZ_i)Z_i,
$$ proving   $RHS=LHS$ and hence confirming (\ref{d-comm-w-flat}).

The form $
\alpha_1$ 
 can now be written as
$$
\alpha_1(Y)=\frac 12(A_g^*(Y)+A(Y))
=A(Y) +\frac 12 (A_g^\ast(Y)-A(Y))
$$
with the second term $\frac 12 (A_g^\ast(Y)-A(Y))$ being exact,
 which
implies that $d^g\alpha_1=d^g A$.
The equation
(\ref{cnd-1}) can now be written as
\begin{equation}
\label{da-1-1}
0=(d^{g} \alpha_1)
(X, Z)(Y)
+(C(Z)Y)X-(C(X)Y)Z
\end{equation}
where $\alpha_1(X)$ is trace-free
and  symmetric with respect to $g$.
This in turn implies that
the map $Y\to C(Z,Y)X-C(X, Y)Z
$ is symmetric,
\begin{equation}
\begin{split}
&\qquad g(C(Z,Y)X-C(X, Y)Z,
W)
\\
& =
g(Y, C(Z,W)X-C(X, W)Z).
\end{split}
\end{equation}
Let $\{Z_i\}$ be an orthonormal basis,
 $Y=Z=Z_i$, and summing over  $i$, we get
%$$
%\tr_g C
%g(X, W)
% -C(X, W)
%=-(n-1)C(X, W)
%$$
%proving that
\begin{equation}
\label{b-c}
(\tr_g C ) g(X, W) + (n-2) C(X, W) =0.
\end{equation}
Taking again the trace we find
\begin{equation}
\label{b-c-1}
\tr_g C  =0.
\end{equation}
Hence $\alpha_0$ satisfies (q2).

Substituting this
into the previous formula we get
$$
(n-2)(C(X, W))=0,
$$
and consequently
\begin{equation}
\label{b-c-2}
C(X, W)=0
\end{equation}
if $n>2 $.

For $n> 2$ it follows
from (\ref{b-c-2}) and (\ref{a-star-a}) that
\begin{equation}
A_g^\ast(Y)=A(Y)
\end{equation}
i.e, $A(Y)$ is symmetric with respect to $g$. Consequently
$$
\alpha_1(Y)=A(Y)=A_g^*(Y).
$$
Hence $\alpha_1(Y)X=A(Y)X=A(X)Y=\alpha_1(X)(Y)$ by Equation (\ref{A-b-sym}).
This proves that $\alpha_1$ satisfies (c1)-(c2), (c4) for $n>2$.
The equation (\ref{da-1-1}) combined with $C=0$ implies further
 $d^g\alpha_1=0$. Hence $\alpha$ satisfies all the conditions to be a cubic form.

%Let $n=2$.
%The symmetric part of
%(\ref{da-1-1}) now reads,
%observing that
% $(d^{g} \alpha_1)(X, Z)$ is
 %symmetric since $\alpha_1$ is symmetric,  that
%$$
%0=(d^{g} \alpha_1)
%(X, Z)
 %+ E(X, Z),$$
%$$ \quad E(X, Z):=
 %\frac 12(X\otimes C(Z) -Z\otimes C(X)
%+ C(Z)^\sharp \otimes X^\flat
%-C(X)^\sharp \otimes Z^\flat).
%$$
%But $ E(X, Z)$, written as matrices
%with $g$ being normalized as the identity matrix, is of the form
%$$
%E(X, Z)=
%\frac 12( (X Z^t-ZX^t) C+ C(ZX^t -XZ^t)).
%$$
%However for $2\times 2$ matrix $C$ with zero trace
%(\ref{b-c-1}),
%this symmetric matrix actually vanishes; indeed, $E(X, Z)$ is of
%the form
%$$
%\begin{pmatrix} 0 &x \\
 %-x  &0\\
%\end{pmatrix}
%\begin{pmatrix} c_1 & c_2\\
 %c_2  &-c_1\\
%\end{pmatrix}
%+
%\begin{pmatrix} c_1 & c_2\\
 %c_2  &-c_1\\
%\end{pmatrix} \begin{pmatrix} 0 &x \\
 %-x  &0\\
%\end{pmatrix}
%=0.
%$$
%Thus $0=d^{g} \alpha_1$. Finally
%taking the trace of (\ref{cnd-1})
 %with respect to the
%second and the third arguments $Z, Y$
%and computing  its inner product with $W$,
%we find, using
%$\tr_g C=0$, that
%\begin{equation*}
%\begin{split}
%0&=\tr_{g, 23}(d^{g} A)(X, W)
%- C(X, W)
%+  \tr_g C g(X, W)
%\\
%&\end{split}
%\end{equation*}
%But $d^g A=0$ so we find $C=0$. Substituting this into (\ref{a-star-a}) we find
%$
%A(Y)^\ast =A(Y)
%$,
%and
%$$
%\alpha_1(Y)=A(Y)=A^*(Y).
%$$
%This proves that  $\alpha_1$ satisfies (c1)-(c4)
%for $n= 2$ also.
\end{proof}

%The following lemma is elementary and is obtained by  elementary computations.
%\begin{lemm+}\label{lemma2.6} Let
 %$$J=\begin{bmatrix}
   %          I_n & 0 \\
     %        0 & -1 \end{bmatrix}$$
%and consider the $SO(n, 1)$-invariant form
%on $E$, $(y, y)_J=g(Y, Y) -\lambda^2$, $y=(Y, \lambda)$.
%Then $(y, y)_J$ is parallel with respect to $\nabla$. In particular
%$J$ induces an involution on
 %$H^1(M, \mathfrak{sl}(n+1,\br))$
%and a decomposition of $\mathfrak{sl}(n+1,\br)=\mathfrak{so}(n,1)\oplus W$, and a decomposition of
%$H^1(M, \mathfrak{sl}(n+1,\br))$: each $\alpha$
%can be written as $\alpha=\alpha_0 +\alpha_1$, with $\alpha_0 \in H^1(M, \mathfrak{so}(n, 1))$ consisting of elements of  the form %$$
%\begin{bmatrix}
%A& B
%\\
%C& 0
%\end{bmatrix}, \quad A=-A^t, B^t=C$$
%and  $\alpha_1$ taking values in
%$W$ consisting of elements of the form
%$$
%\begin{bmatrix}
%A& B
%\\
%C& D
%\end{bmatrix}, \quad A=A^t, B^t=-C, D=-\tr A.$$
%Here the transpose $B^t$ is computed with respect to $g$.
%\end{lemm+}

We prove now that the map from $\alpha$ to the cubic form
in  Lemma \ref{lemma2.5} is injective.
\begin{lemm+}\label{2.7}Let $n>2$ and suppose $\alpha\in
H^1(\pi_1(M),\rho, \mathfrak{sl}(n+1,\br))$.
Then the map
$\alpha\mapsto \alpha_1$ from $H^1(\pi_1(M),\rho, \mathfrak{sl}(n+1,\br))$
to the cubic form $\alpha_1$ is injective.
\end{lemm+}
\begin{proof}
In Lemma \ref{lemma2.5}, we showed that if $\alpha$ is represented by a $\mathfrak{sl}(E)$-valued one form with $B=D=0$, then
$C=0$ and the associated cubic form is $\alpha_1(Y)=A(Y)$. Hence if the associated cubic form vanishes, $A=0$. This implies that
$\alpha$ is represented by an exact 1-form by Equation (\ref{newA}), hence it is  a zero element in the cohomology.
\end{proof}
This finishes the proof of Theorem \ref{2.2}.

Now we  prove the
 Weil's local rigidity theorem,
 $ H^1(\Gamma, \mathfrak{so}(n, 1))=0$, $n>2$,
\cite[Chapter VII]{Rag-book}, 
\cite{Weil}  as an application of our technique.
For $n=2$, it is well-known that the cohomology  $ H^1(\Gamma, \mathfrak{so}(2, 1))$
is   determined by quadratic forms satisfying (q1-q3),
namely real part of holomorphic quadratic forms; see e.g.
\cite{Eichler}.
\begin{theo+}\label{Weil} Let $M=\Gamma\backslash SO^0(n, 1)/SO(n)$
be a compact hyperbolic manifold.
\begin{enumerate}
\item
If $n>2$ then   $ H^1(\Gamma, \mathfrak{so}(n, 1))=0$.
\item If $n=2$ then  $ H^1(\Gamma, \mathfrak{so}(n, 1))$ is given
by the space of quadratic forms satisfying (q1-q3).
\end{enumerate}
\end{theo+}
\begin{proof}  Let $\alpha $
represent an element in $ H^1(\Gamma, \mathfrak{so}(n, 1))$ viewed as an element in $\Omega^1(M,\mathfrak{so}(n, 1))$. The elements 
in $\mathfrak{so}(n, 1)$ are of the form
$\begin{bmatrix}
a& b\\
c&0
\end
{bmatrix}$ with $a=-a^*$ and $b=c^T$ with respect
to the Euclidean product in $\mathbb R^n$ 
as a subspace of the Lorentz space $\mathbb R^{n+1}$. 
The $1$-form $\alpha$ takes then the form
$$
\alpha =
\begin{bmatrix}
A& B\\
C&0
\end{bmatrix}
$$ 
with $g(A(X)Y, Z)=-g(Y, A(X)Z)$,
$g(B(X), Y)= C(X)(Y)=C(X, Y)$, where $g$ is the given
hyperbolic metrix on $M$.
The (1,1)-entry $A$ of $\alpha$
is skew-symmetric.
Now from Equations (\ref{ext-diff}), (\ref{dA}) and the fact that $\nabla^T=\nabla^g$ for hyperbolic manifold,   $(1, 1)$-entry of the condition $d^g \alpha=0$
is
\begin{equation}
\label{1-1-d}
(d^g A)(X, Z) + X\otimes C(Z) -Z\otimes C(X) + B(X) \otimes Z^\sharp
-B(Z)\otimes X^\sharp=0
\end{equation}
as two form acting on $(X, Z)$. Here $X\otimes Z^\sharp $
is the rank-one map $Y\mapsto g(Y, Z) X$.
We shall also need some Hodge theory. Equip $so(n, 1)$
with the $SO(n)$-invariant positive inner product induced
from the standard Euclidean inner product in $\mathbb R^{n+1}$, $(y,
y)_E=\Vert
Y\Vert^2_g +\lambda^2$, $y=(Y, \lambda)$.
Take a harmonic  one form representing $\alpha$.
Then the cohomology class $\alpha$ satisfies also
the coboundary condition $\nabla^\ast \alpha=0$.
To write
down the formula for $\nabla^\ast $
we observe that
$$
\nabla_X=
\begin{bmatrix}
\nabla^g_X &0
\\
0&L_X
\end{bmatrix}
+
\begin{bmatrix}
0&X
\\
X^\sharp &0
\end{bmatrix}
$$
is a sum of two terms, the first
preserving the
Euclidean inner product $(y, y)_E$,
whose adjoint can be found
by standard formula
(see e.g. \cite[p.2]{Matsu}), whereas
the second part is self-adjoint.
Thus
$-\nabla^\ast \alpha$ is given by
$$
\sum_j (\delta_{X_j} \alpha) (X_j)
$$
where
$$
\delta_X
=\begin{pmatrix} \nabla^g_X & -X\\
-X^\sharp & L_X
\end{pmatrix}.
$$
More precisely, $(\delta_X\alpha)(Z)$ is given, 
for any testing section $y=(Y, \lambda)$,  by
the Leibniz rule
\begin{equation*}
  \begin{split}
&\qquad (\delta_X\alpha)
(Z)y\\
&=\begin{pmatrix} \nabla^g_X & -X\\
-X^\sharp & L_X
\end{pmatrix}
(\alpha(Z)y)- \alpha(Z) \left(\begin{pmatrix} \nabla^g_X & -X
\\
-X^\sharp & L_X
\end{pmatrix} y\right)-\alpha(\nabla_X^g Z)y.
  \end{split}
\end{equation*}
(The sum
$\sum_j (\delta_{X_j} \alpha) (X_j)$ is well-defined but
not the individual terms.)
When acting on the section $y=(0, 1)$ we find
\begin{equation}
  \label{dual}
\sum_j( A(X_j)X_j +(\nabla^g_{X_j} B)(X_j))=0.
\end{equation}

It now follows from Equation (\ref{newA}), Lemmas 2.3 and 2.5 (keeping track of the change
of forms) that
$A_1(X):=A(X) +\nabla_X^g B + (Xb)I +f X$
is symmetric and trace-free and satisfies the condition (c1-c4).
Note
here while performing computations as in Lemmas 2.3-2.5
we use forms $u$ with values in $\mathfrak {sl}(n+1)$
instead of in  $\mathfrak {so}(n, 1)$, however all
we need is that $d^\nabla d^\nabla =0$, i.e. we will show that $\alpha$ vanishes identically.
The trace free condition  and $$A(X)Y +\nabla_X^g B(Y)+(Xb)Y +f(Y) X=A(Y)X +\nabla^g_Y B(X)+(Yb)X+ f(X)Y$$ imply  that
the map $Z\to A_1(Z)Y=A(Z)Y +(\nabla_Z B)(Y) +(Zb)Y + f(Y)Z$
is trace free.
The symmetric relation implies
$$
g(A(Z)Y+ (\nabla_Z B)(Y) +(Zb)Y + f(Y)Z, W)$$ $$
=g(A(Z)W+ (\nabla_Z B)(W) + (Zb)W + f(W)Z, Y).
$$
We take $\{Z_j\}$ a local orthonormal frame  and
put $Z=W=Z_j$ in the above equation. Summing over $j$
we find, in view of
(\ref{dual}) that the right hand side is
$$
RHS=\sum_j g(A(Z_j)Z_j+ (\nabla_{Z_j} B)(Z_j) +(Z_jb)Z_j + f(Z_j)Z_j, Y)
=(Yb) +f(Y)
$$
and
$$
LHS=\tr (A_1(\cdot)Y)=0.
$$
Namely the one-form $Yb+f(Y)=0$. But $f(Y)=(db)(Y)=Yb$
by (\ref{newA}),
so $0=2f(Y)$,  and $db= f=0$.
This implies  in turn that
$$A_1(X)=
A(X) +\nabla_X^g B + (Xb)I +f X
=A(X)+\nabla_X^g B$$
is symmetric.
We write $B=B^0+B^1$, the symmetric and respectively
the skew symmetric part of $B$.
Since $A$ is skew symmetric,  the skew symmetric part of
$A_1(X)$ must vanish, that is 
\begin{equation}
\label{a+nabla-b}
A(X)+\nabla_X^gB^1=0.
\end{equation}
This implies
in turn $A(X)=-\nabla_XB^1$ is exact.
Thus $d^g A= 0$, and the relation (\ref{1-1-d}) becomes
$$
X\otimes C(Z) -Z\otimes C(X) + B(X) \otimes Z^\sharp
-B(Z)\otimes X^\sharp=0.
$$
Let $\{Y_i\}$ be an orthonormal basis, put $Z=Y_i$ and let the above
act on $Y_j$. Taking the sum
and using  $B^t=C$ we find  as in the proof of Lemma 2.5,
that
$$ \sum_i g(X\otimes C(Y_i,Y_i)-Y_i\otimes C(X,Y_i)+B(X)g(Y_i,Y_i)-B(Y_i)g(X,Y_i), W)=0$$ i.e.
$$
\tr_g C g(X, W) - C(X, W) + n C(X, W) - C(X, W)=0.
$$
But the same proof above implies that
\begin{equation}
\label{crit-eq} 
\tr_gC=0, \,\, (n-2)C(X, W)=0.
\end{equation}

Now let $n>2$. Thus $C(X, W)=0$. Then $B=C^t=0$  and using  the symmetric condition on $A_1=A$
we find that $A$ is symmetric and thus $A=0$.  This proves (1).

Let $n=2$.  We consider the $\mathfrak{so(2, 1)}$-valued section
$$
u=
\begin{pmatrix}
B^1&0
\\
0&0
\end{pmatrix}.
$$
Using the formulas (\ref{nabla-u}) 
and
(\ref{a+nabla-b})
we find 
$$
\alpha +\nabla u
:X\mapsto \alpha(X)
+\nabla_X u
=\begin{bmatrix}
0& B^0
\\
C^0&0
\end{bmatrix}
$$
where $C^0$ is the symmetric part of $C$ and $(B^0)^t=C^0$. So replacing
$\alpha$ by
$\alpha+\nabla u $ we may assume that $A=0$,  $B$ is symmetric and $B^t=C$.
A direct calculation using (\ref{da=0}) and $y=(Y,0)$  gives
\begin{equation*}
\begin{split}
&0=d^\nabla \alpha(X,Z)Y\\
&=\begin{bmatrix}
(C(Z)Y)X-g(X,Y)B(Z)- (C(X)Y)Z+ g(Z,Y)B(X) \\
X(C(Z)Y)-C(Z)\nabla^g_XY- Z(C(X)Y)+ C(X)\nabla^g_ZY - C([X, Z])Y \end{bmatrix}.
\end{split}
\end{equation*}
But by the formula (\ref{dc})
$$X(C(Z)Y)-C(Z)\nabla^g_XY- Z(C(X)Y)+ C(X)\nabla^g_ZY - C([X, Z])Y=(d^gC)(X,Z)Y$$

Hence  $d^\nabla \alpha=0$ gives $d^g C=d^g B=0$. 
We have thus that $B$ is symmetric, $\tr_gB=\tr_g C=0$
by (\ref{crit-eq}) and $d^g B=0$, namely $B$ satisfies 
(q1-q3).

This completes the proof.
\end{proof}

Recall \cite{Hitchin} that
the Hitchin component,  denoted by  $\chi_H(\pi_1(S),
SL(3,\br))$,
is a connected component in the character variety
 $\mathrm{Hom}(\pi_1(S), SL(3, \mathbb R))//SL(3, \mathbb R)$
containing the realization of $\pi_1(S)$
as a subgroup of $SL(2, \mathbb R)$
composed with the irreducible
representation of $SL(2, \mathbb R)$
on $\mathbb R^3$. The tangent space of  $\chi_H(\pi_1(S),
SL(3,\br))$ at
these specific Fuchsian points can
 be obtained from the general
theory in \cite{Hitchin}. See also \cite{KZ}.


\begin{thebibliography}{99}

%\bibitem{Ahlfors-1961}L. V. Ahlfors, \emph{Some remarks on Teichm\"u{}ller's space of Riemann surfaces,} Ann. of Math. 74 (2) (1961),  171-191.

%\bibitem{Ahlfors-1961/62} L. V. Ahlfors, \emph{Curvature properties of Teichm\"u{}ller's space,} J. Analyse Math. 9 (1961/1962),  161-176.

\bibitem{Be1}Y. Benoist, \emph{Automorphismes des c\^ones convexes,} Invent. Math. 141 (1) (2000), 149-193.

\bibitem{Be}Y. Benoist, \emph{Convexes divisibles III,}  Ann. Sci. de l'ENS 38 (2005), 793-832.

%\bibitem{Bob} {B. Berndtsson}, \emph{Curvature of vector bundles associated to holomorphic fibrations}, {Ann. of Math. (2)}, \textbf {169} ({2009} {531--560}.


%\bibitem{Bo-mz} {B. Berndtsson}, \emph{Strict and nonstrict positivity of direct image bundles}, {Math. Z.},  \textbf {269} ({2011}), {1201--1218}.

%\bibitem
%{BCLS} M. Bridgeman, R. Canary, F. Labourie, and A. Sambarino, \emph{The pressure metric for convex representations,} GAFA. \textbf{2} (2015), 1089-1179.

\bibitem{CG} S. Choi and W. Goldman, \emph{Convex real projective structures on closed surfaces are closed,} Proc. A.M.S. 118 (1993), 657-661.

%\bibitem{Co}K. Corlette, \emph{Flat G-bundles with canonical metrics}, J.D.G. 28 (1988), no. 3, 361-382.

\bibitem{Eichler}
M. Eichler,
\emph{
Eine Verallgemeinerung der Abelschen Integrale,}
Math. Z. 67 (1957), 267-298.




\bibitem{Hitchin}
N. Hitchin,
\emph{Lie groups and {T}eichm\"uller space},
{Topology},
\textbf{31}(1992),
 {449--473}.

\bibitem{Kim}I. Kim, \emph{Rigidity and deformation spaces of strictly convex real projective structures on compact manifolds},  J.D.G. 58 (2001), 189-218.
    Erratum,  J.D.G. 86 (2010), 189.

\bibitem{KZ} I. Kim and  G. Zhang, \emph{Eichler-Shimura isomorphism for complex hyperbolic lattices}, submitted.

\bibitem{KZ-hitchin}  I. Kim and  G. Zhang,
\emph{K\"ahler metric on the space of convex real projective structures on surface}, J.D.G., to appear. 


\bibitem{La}
F. Labourie, \emph{Flat projective structures on surfaces and cubic holomorphic differentials}, Pure Appl. Math. Q, \textbf{3} (2007), 1057-1099.






\bibitem{Loftin}
 J. Loftin, \emph{Affine spheres and convex $RP^n$ manifolds},
American Journal of Math.
\textbf{
123}
 (2) (2001), 255-274.


\bibitem{Matsu} Y. Matsushima,
\emph{Vector bundle valued
harmonic forms and immersions of
Riemannian manifolds},
Osaka J. Math.,
\textbf{
8} (1971), 1-13.

\bibitem{Rag-book}
M.S. Raghunathan, \emph{Discrete subgroups of {L}ie groups}, Springer-Verlag, New York,
  1972, Ergebnisse der Mathematik und ihrer Grenzgebiete, Band 68. \MR{0507234
  (58 \#22394a)}

\bibitem{Weil}A. Weil, \emph{Discrete subgroups of Lie groups, II}, Ann. of
Math \textbf{75} (1962), 97-123.

\end{thebibliography}
\end{document}